\documentclass{article}
\usepackage{fullpage}
\usepackage{amsmath}   
\usepackage{amssymb}
\usepackage{amsthm}
\usepackage{ulem}
\usepackage{epsfig}

\title{A space-time integral estimate for a large data semi-linear wave equation on the Schwarzschild manifold. }
\author{P. Blue \and A. Soffer}

\newtheorem{theorem}{Theorem}
\newtheorem{lemma}[theorem]{Lemma}

\newcommand{\Reals}{\mathbb{R}}
\newcommand{\man}{\mathcal{M}}

\newcommand{\dt}{\partial_t}
\newcommand{\dr}{\partial_\rho}
\newcommand{\slap}{\Delta_{S^2}}
\newcommand{\sgrad}{\nabla_{S^2}}

\newcommand{\dthree}{d\rho d^2\omega}

\newcommand{\chiHardy}{\chi}
\newcommand{\chispec}{\chi_{1}}
\newcommand{\epsilonspec}{\epsilon_{1}}

\newcommand{\confchg}{E_\mathcal{C}}
\newcommand{\intTslice}{\int_{ \{t\}\times\man} }
\newcommand{\intTslab}{\int_{ (t_1,t_2)\times\man} }

\newcommand{\chiLight}{ {{\tilde{\chi}}_\text{lc}} }
\newcommand{\chiLightn}[1]{ {{\tilde{\chi}}_{\text{lc}, #1}} }
\newcommand{\Nlightcone}{ {N_\text{lc}}}
\newcommand{\Clightcone}{ {C_\text{lc}}}

\newcommand{\underbracetag}[2]{\underbrace{#1}_{#2}}

\newcommand{\Op}{A}

\newcommand{\GEFA}{\psi_0}
\newcommand{\GEF}{\phi_0}
\newcommand{\LEF}{\phi_\text{GS}}
\newcommand{\LEV}{E_{\text{GS}}}

\newcommand{\dx}{\partial_x}
\newcommand{\dthreex}{dxd^2\omega}

\newcommand{\Hide}[1]{}

\begin{document}

\maketitle

\begin{abstract}
We consider the wave equation $(-\dt^2+\dr^2 -V -V_L(-\Delta_{S^2})) u = fF'(|u|^2) u$ with $(t,\rho,\theta,\phi)$ in $\Reals\times\Reals\times S^2$. The wave equation on a spherically symmetric manifold with a single closed geodesic surface or on the exterior of the Schwarzschild manifold can be reduced to this form. Using a smoothed Morawetz estimate which does not require a spherical harmonic decomposition, we show that there is decay in $L^2_{\text{loc}}$ for initial data in the energy class, even if the initial data is large. This requires certain conditions on the potentials $V$, $V_L$ and $f$. We show that a key condition on the weight in the smoothed Morawetz estimate can be reduced to an ODE condition, which is verified numerically. 
\end{abstract}

We consider the following defocusing wave equation:
\begin{align}
 ( -\dt^2 +\dr^2 -V -V_L(-\slap) ) u =& f(\rho) F'(|u|^2) u \label{eNLW}\\
u(0)=& u_0 \nonumber \\
\dt u(0)=& u_1 \nonumber
\end{align}
for a function $u:\Reals\times\man\rightarrow \Reals$ with $(t,\rho,\theta,\phi)\in\Reals\times\man=\Reals\times\Reals\times S^2$. Our goal is to show that, under conditions \ref{cRadial}-\ref{cSpectralCondition} on the potentials (see below), 
\begin{align*}
\| u\|_{L^2_{\text{loc}}} \rightarrow 0
\end{align*}
if the initial data $(u_0,u_1)$ has finite (but possibly large) energy. In fact, we show the stronger result that $\int_0^t \int_\man (1+\rho^2)^{-2} |u|^2 dt\dthree < C E[u_0,u_1]$. 

Equation \eqref{eNLW} describes several interesting geometric wave equations. On a spherically-symmetric, $3$-dimensional, warped-product of $\Reals$ and $S^2$, Riemannian manifold, the metric is $ds^2 =d\rho^2 + r(\rho)^2 d\omega^2$, and if $r$ has a unique minimum, it corresponds to a closed geodesic surface. This is the first problem we consider, and, by translation, we may assume that the closed geodesic surface occurs at $\rho=0$. The semi-linear wave equation $(-\dt^2+\Delta_{g})\tilde{u}=|\tilde{u}|^{p-1}\tilde{u}$ can be reduced to \eqref{eNLW} by setting $u=r\tilde{u}$, in which case $V=r''/r$, $V_L=r^{-2}$, and $f=r^{1-p}=V_L^{(p-1)/2}$. The conditions \ref{cRadial}-\ref{cSpectralCondition} on the potentials are, therefore, conditions on $r$. These conditions are not vacuous, since, for example $r=1+\rho^2$ generates a wave equation which satisfies these conditions. Similarly the wave equation on the exterior region of the $3+1$-dimensional, Schwarzschild manifold also can be reduced to \eqref{eNLW} if $p=3$, and this is the second case we consider. Both cases are explained in more detail in \cite{BlueSoffer4}. 

The purpose of this note is to show that there is decay even for large data (defocusing) semi-linear problems. 
In \cite{BlueSoffer4}, we considered the linear equation (\eqref{eNLW} in which the right-hand side is absent) in both the Riemannian and Schwarzschild cases, and the non-linear problem with small initial data in the Riemannian case. In these cases, we were able to show that the space-time integral $\|\tilde{u}\|_{L^4(dt\dthree)}$ is controlled by weighted $H^{1+\epsilon}$ norms. This built on previous work \cite{BlueSterbenz}, in which the $L^\infty$ norm was shown decay if weighted $H^{3}$ norms were bounded (and the initial data was small in the non-linear case). Similar results, with better decay estimates at the event horizon (an important part of the Schwarzschild manifold, corresponding to $\rho\rightarrow-\infty$) have also been proven \cite{DafermosRodnianski}. 

We make the following assumptions on the potentials and non-linear terms: 
\begin{enumerate}
\item {}[Radial potentials] $V$ and $V_L$ are smooth functions of $\rho$ only. \label{cRadial}
\item {}[Positive energy] $V$ and $V_L$ are positive. \label{cPositiveEnergy}
\item {}[Unique, unstable, trapped geodesic surface] $V_L$ (which we refer to as the angular potential) has a unique critical point, which is a maximum, and occurs at $\rho=0$. We refer to $\rho=0$ as the trapped ray, trapped geodesic surface, or the orbiting null-geodesic surface. \label{cVLuniquepeak}
\item {}[Derivative of $V$ well-behaved] $\dr V$ is bounded and goes to zero at $\pm\infty$. \label{cVWellBehaved}
\item {}[Semi-linearity] There's a $p\in (1,3]$ such that $F(x)=\frac{|x|^{\frac{p+1}{2}} }{p+1}$. 
\item {}[Compatible semi-linear coefficient and angular potential] $f=V_L^{\frac{p-1}{2}}$. \label{cCompatibleSemiLin}
\end{enumerate}
Throughout this note, we use $C$ to denote an arbitrary constant which may change from line to line, $\langle\cdot,\cdot\rangle$ as the $L^2$ inner product (with respect to $\dthreex$ unless otherwise stated), and $\|\cdot\|$ as the $L^2$ norm. 

\section{Analysis using the method of multipliers}

It is well-known that there is a conserved energy, as in $\Reals^{n+1}$, 
\begin{align}
E[u,\dt u]
=& \int_{ \{t\}\times\man} |\dt u|^2 + |\dr u|^2 + V|u|^2 + V_L |\sgrad u|^2 + 2 f(\rho) F(|u|^2) \dthree . 
\end{align}

Following our earlier work \cite{\Hide{LabaSoffer,}BlueSoffer4}, we introduce a radial derivative operator, pointing away from the trapped ray
\begin{align*}
\gamma =    g\dr + g'/2 
\end{align*}
in terms of a $C^3$ radial function $g=g(\rho)$ satisfying 
\begin{enumerate}
\setcounter{enumi}{6}
\item {}[Increasing] $g'(\rho)=\dr g>0$, \label{cgIncreasing}
\item {}[Centered about the trapped ray] $g(0)=0$,  \label{cgAwayFromZero}
\item {}[Energy bounded] $g$ is bounded, \label{cgBounded}
\item {}[Inverse polynomial behaviour] there's a power $\alpha\in[-2,-1)$, and for $i\in\{1,2,3\}$, there are constants $c_i$ and $C_i$ such that (for $|\rho|$ sufficiently large) $c_i \rho^{\alpha+1-i}\leq  g^{[i]}\leq C_i \rho^{\alpha+1-i}$, \label{cgPolynomial} and 
\item {}[Positive spectral condition] \label{cSpectralCondition}There's an $\epsilonspec>0$ and a non-negative function $\chispec$ which is positive in a non-empty open set, such that $\forall \psi\in C^\infty(\{t\}\times\man):$
\begin{align*}
\intTslice 2g'|\dr\psi|^2 + (-g'''/2 -gV')|\psi|^2 \dthree 
\geq \epsilonspec \intTslice g'|\dr\psi|^2 + \chispec |\psi|^2 \dthree .
\end{align*}
This is equivalent to saying that the (self-adjoint) operator $\Op = -(2-\epsilonspec)\dr g'\dr -g'''/2 -gV' -\epsilonspec\chispec$ has non-negative spectrum. 
\end{enumerate}
Because condition \ref{cSpectralCondition} involves $V$, the existence of $g$ satisfying conditions \ref{cgIncreasing}-\ref{cSpectralCondition} is a condition on $V$. In the Riemannian case, with $r=1+\rho^2$, since $V=2/r^2$ has a unique maximum at $\rho=0$, the arguments in \cite{BlueSterbenz} or \cite{BlueSoffer4} shows the existence of a $g$ satisfying conditions \ref{cgIncreasing}-\ref{cSpectralCondition} (technically, the functions in \cite{BlueSterbenz} and \cite{BlueSoffer4} have $g'''$ as a distribution, not a continuous function, but this is sufficient for the analysis). 

Applying the method of multipliers to $\gamma u$, we have
\begin{align*}
(\langle \dot{u},\gamma u\rangle+ \langle \gamma u,\dot{u}\rangle) |_{t_1}^{t_2}
=& \intTslab 
\underbracetag{2g'|\dr u|^2}{\text{Term $I$}} 
+\underbracetag{(-g'''/2 -gV')|u|^2}{\text{Term $II$}} 
+\underbracetag{ (-gV_L')|\sgrad u|^2 }{\text{Term $III$}} \dthree dt\\
& + \underbracetag{\int_{t_1}^{t_2} \left(\langle\gamma u,  f F'(|u|^2) u\rangle + \langle f F'(|u|^2) u,\gamma u\rangle \right) dt}{\text{Term $IV$}} .
\end{align*}

We now show that the integrand on the right is positive. By condition \ref{cSpectralCondition}, terms $\text{I}$ and $\text{II}$ are positive when taken together, and, in fact, dominate $(1+\rho^2)^{-2}|u|^2$ through the Hardy estimate in lemma \ref{lHardy}. The contribution from term $\text{III}$ is positive because by conditions \ref{cVLuniquepeak} and \ref{cgAwayFromZero}, $g$ changes from negative to positive at the same point $-V_L'$ does. At this point, we analyse the structure of the non-linear contribution, term $\text{IV}$ in greater detail than in \cite{BlueSoffer4}. 
\begin{align*}
\langle\gamma u, & f F'(|u|^2) u\rangle + \langle f F'(|u|^2) u,\gamma u\rangle \\
=& -\intTslice g \dr( f F'(|u|^2) ) |u|^2 \dthree \\
=& -\intTslice g \dr( f |u|^{p-1} ) |u|^2 \dthree \\
=& -\intTslice g \frac{p-1}{p+1} f^{\frac{-2}{p+1}} \dr (f^{\frac{p+1}{p-1}} |u|^{p+1}) \dthree \\
=&  \frac{p-1}{p+1} \intTslice f^{\frac{p+1}{p-1}} |u|^{p+1} \dr (g f^{\frac{-2}{p-1}}) \dthree\\
=& (p-1) \intTslice (f (\dr g) -\frac{2}{p-1} g (\dr f)) F(|u|^2)\dthree \\
=& (p-1) \intTslice \left( V_L^{\frac{p-1}{2}} (\dr g) - g  V_L^{\frac{p-3}{2}}(\dr V_L) \right) F(|u|^2)\dthree .
\end{align*}
Since $V_L$ is positive and $\dr g>0$, the first term in brackets is positive. Since $g$ changes sign at the maximum of $V_L$, $-g\dr V_L >0$ is also positive. Thus the integrand is positive. 

We now turn to the left-hand side. Since $g$ is bounded, $\|g\dr u\|^2< CE$. Since $g'<(1+|\rho|)^{-1}$, by the Hardy estimate, $\|(\dr g)u\|^2<C E$. Thus, the left hand-side is bounded by the energy, and
\begin{align}
CE 
>& \intTslab \frac{|u|^2}{(1+\rho^2)^2} dt\dthree \label{eMainResult}\\
&+\intTslab g(\dr V_L) |\sgrad u|^2 + \left( V_L^{\frac{p-1}{2}} (\dr g) - g  V_L^{\frac{p-3}{2}}(\dr V_L) \right) F(|u|^2)  dt\dthree .\nonumber
\end{align}
Taking $t_1=0$ and $t_2\rightarrow\infty$ gives the desired result. 

We now prove a smooth Hardy estimate. Essentially, since $\man=\Reals\times S^2$ is a $3$-dimensional manifold, one expects the Hardy estimate from $\Reals^{3}$, $\|u/|x|\| \leq \| \nabla u\|$ to hold; however, because there is no origin for $\Reals\times S^2$, a little extra control is required. This result is taken from \cite{BlueSterbenz}. 

\begin{lemma}[Smooth Hardy]
\label{lHardy}
For any $\alpha\geq0$ and non-negative function $\chiHardy$ which is positive in a non-empty, open set, there's a constant $C$, such that
\begin{align*}
\int \frac{|\dr u|^2}{ (1+\rho^2)^{\frac{\alpha}{2}} } +\chiHardy |u|^2 \dthree 
> C \int \frac{1}{(1+\rho^2)^\frac{\alpha+2}{2}} |u|^2 \dthree 
\end{align*}
In particular, when $\alpha=0$, since $V$ is a smooth positive function, $E>C\int (1+\rho^2)^{-1}|u|^2 \dthree$. 
\begin{proof}
For $\rho_1>0$, 
\begin{align*}
\frac{|u(\rho_1)|^2}{(1+\rho_1)^{\alpha+1}} - |u(0)|^2 
=& \int_0^{\rho_1} \dr \frac{|u|^2}{(1+\rho)^{\alpha+1}} d\rho \\
=&\int_0^{\rho_1} \frac{2u\dr u}{(1+\rho)^{\alpha+1} } -(\alpha+1)\frac{|u|^2}{(1+\rho)^{\alpha+2} } d\rho\\
\leq& \frac{\alpha+1}{2}\int_0^{\rho_1} \frac{|u|^2}{(1+\rho)^{\alpha+2}} d\rho 
+ \frac{2}{\alpha+1} \int_0^{\rho_1}\frac{|\dr u|^2}{(1+\rho)^{\alpha}} d\rho 
-(\alpha+1) \int_0^{\rho_1} \frac{|u|^2}{(1+\rho)^{\alpha+2}} d\rho\\
\int_0^{\rho_1} \frac{|u|^2}{(1+\rho)^{\alpha+2}} d\rho
\leq& \frac{4}{(\alpha+1)^2} \int_0^{\rho_1}\frac{|\dr u|^2}{(1+\rho)^{\alpha}} d\rho 
+\frac{2}{\alpha+1} |u(0)|^2 .
\end{align*}
We take the limit $\rho_1\rightarrow\infty$. Since (for any exponent $\beta\geq0$) $(1+\rho)^{\beta}$ is equivalent to $(1+\rho^2)^{\beta/2}$ on $[0,\infty)$, the powers of $(1+\rho)$ can be replaced by $(1+\rho^2)^{1/2}$. By symmetry, the same result holds on $(-\infty,0]$. Since $(1+\rho^2)^{-\beta}$ is uniformly equivalent to $(1+(\rho-\rho_0)^2)^{-\beta}$ for $\rho_0$ in a finite interval, the $|u(0)|^2$ term can be replaced by $|u(\rho_0)|^2$ in any fixed interval. By integrating the estimate over $\rho_0$ with $\rho_0$ in a bounded open set in which $\chiHardy$ is positive, the desired result holds. 
\end{proof}
\end{lemma}

\section{Using numerical ODE solutions to verify the spectral condition
}

In this section, we (i) show that the spectral condition can be reduced to showing that there are no zeroes for the solution to an associated ODE and (ii) numerically verify this condition for some choice of function $g$ on the Schwarzschild manifold. 

\subsection{Reducing the spectral condition to an ODE condition }

\begin{lemma}
The spectral condition, \ref{cSpectralCondition}, follows from the following conditions: 
\begin{enumerate}
\setcounter{enumi}{11}
\item
{}[ODE condition] There is a smooth, non-negative solution $\GEFA$ to \label{cODE}
\begin{align}
(-(2-\epsilonspec)\dr g'\dr -g'''/2 -gV'-\epsilonspec\chispec)\GEFA = 0 .
\label{ecODE}
\end{align}
\end{enumerate}
\begin{proof}
We introduce several new variables to make this problem simpler. First, we change variables to $x=x(\rho)$ defined implicitly by
\begin{align*}
\frac{\text{d} \rho}{\text{d}x} =& g' ,&
x(0)=&0 .
\end{align*}
Since $g'$ is positive, $x$ is well defined. By condition \ref{cgPolynomial}, $x$ grows polynomially in $\rho$. Let 
\begin{align*}
\Op=& -(2-\epsilonspec)\dr g'\dr -g'''/2 -gV' -\epsilonspec\chispec\\
B=& -(2-\epsilonspec)\dx^2 +W \\
W=& g'(\rho(x)) \left(-\frac{g'''(\rho(x))}{2} -g(\rho(x))V'(\rho(x)) -\epsilonspec\chispec(\rho(x)) \right) .
\end{align*}
Note that, by conditions \ref{cVWellBehaved}, \ref{cgBounded}, and \ref{cgPolynomial}, $W\rightarrow 0$ as $x\rightarrow\pm\infty$. If $\psi\in C^\infty(\rho,\omega)$, then for $\phi\in C^\infty(x,\omega)$ given by $\phi(x,\omega)=\psi(\rho(x),\omega)$, 
\begin{align*}
\Op\psi
=& \frac{1}{g'} B\phi , \\
\int \bar{\psi} \Op \psi \dthree
=& \int \bar{\phi} B \phi \dthreex .
\end{align*}
Thus, the spectrum of $\Op$ and $B$ are the same, and if $\Op \psi =0$ with $\psi$ nowhere zero, then $B\phi=0$ with $\phi$ nowhere zero. Let $\GEF$ be a non-negative solution to $B\GEF=0$. 

If $\phi\in C^\infty_0$, then, since $\GEF$ is positive, we can write it as $\phi=\GEF u$ with $u\in C^\infty_0$. Consider the expectation value of $B$ with respect to $\phi$. In this paragraph, we use $c=-(2-\epsilonspec)<0$. 
\begin{align*}
\langle \phi,B\phi\rangle
=& \int \GEF\bar{u} (c\dx^2 +W) \GEF u \dthreex \\
=& \int \GEF\bar{u} \left( (c\dx^2 \GEF)u +2c(\dx\GEF)(\dx u) +\GEF(c\dx^2 u) + W\GEF u \right) \dthreex \\
=& \int \GEF\bar{u} \left( (B\GEF)u +2c(\dx\GEF)(\dx u) +\GEF(c\dx^2u) \right) \dthreex \\
\int \GEF\bar{u} \GEF (c\dx^2 u) \dthreex
=& \int \GEF^2 \left( \dx (\bar{u}(c\dx u)) -(\dx\bar{u}c \dx u) \right)\dthreex \\
=&\int -2\GEF(\dx\GEF)\bar{u}(c\dx u) -\GEF^2(\dx\bar{u}c \dx u) \dthreex \\
\langle \phi,B\phi\rangle 
=& \int -c\GEF^2 |\dx u|^2 \dthreex \\
\geq& 0
\end{align*}
This proves that $B$ and, hence, $A$ have non-negative spectrum.

\Hide{We now consider the negative part of the spectrum of $B$ using results on the spectrum and ground state of Schrodinger operators \cite{LiebLoss}. Since $B$ is a Schrodinger operator, the negative spectrum is discrete. Furthermore, since $W$ is bounded below, if there is any negative spectrum, there is a ``ground state'' energy, $\LEV$ which is the minimum $E$ for which $B\phi= E\phi$ has a solution with $\phi\in L^2$. The function $\LEF$ associated with $\LEV$ is unique and nonnegative. Furthermore, since $W$ is smooth and decaying, $\LEF$ and $\GEF$ are smooth, and $\LEF$ and its derivatives decay exponentially at $\pm\infty$. 

Since $\GEF$ is polynomially bounded, we can consider $\int \GEF\LEF \dthreex$. We may also integrate by parts since the decay for $\LEF$ and its derivatives dominate the possible growth of $\GEF$ and its derivatives. 
\begin{align*}
\LEV \int \GEF\LEF \dthreex
=&\int \GEF(B\LEF) \dthreex \\
=&\int \GEF(-2\dx^2+W)\LEF \dthreex \\
=&\int ((-2\dx^2+W)\GEF) \LEF \dthreex \\
=&0 .
\end{align*}
Since $\LEV<0$ by assumption, $\LEF$ is a non-negative function, and $\GEF$ is a positive functions, this is a contradiction, and there can be no negative spectrum. Thus $B$ and hence $\Op$ are positive operators. }
\end{proof}
\end{lemma}

\subsection{Numerical verification of the ODE condition on the Schwarzschild manifold}

To verify condition \ref{cODE}, one can solve the ODE numerically. There are two potential problems with this method: first, the accuracy of the numerical solution, and second, solutions can only be found on a finite range. Because of the power of computers, we will expect that numerical solutions are accurate and ignore the first problem. We demonstrate how to control the errors caused by the second problem by considering the wave equation on the Schwarzschild manifold. 

On the Schwarzschild manifold, there is another radial coordinate $r\in(2M,\infty)$ defined implicitly, and used to define the potentials, 
\begin{align*}
\frac{\text{d}r}{\text{d}\rho} =& (1-\frac{2M}{r}) ,&
r(0)=&3M, &
V=&\frac{2M}{r^3}(1-\frac{2M}{r}), &
V_L=& \frac{1}{r^2}(1-\frac{2M}{r}) .
\end{align*}
From the definition of $r$, it follows that $1-2M/r$ decays exponentially in $\rho$ as $\rho\rightarrow-\infty$ and $r/\rho\rightarrow 1$ as $\rho\rightarrow\infty$. Let 
\begin{align*}
g=&\int_0^\rho \frac{1}{1+b\tau^2} d\tau , 
\end{align*}
The function $g'$ decays like $b^{-1}\rho^{-2}$, $g'''$ decays like $6b^{-1}\rho^{-4}$, and $gV'$ decays exponentially as $\rho\rightarrow-\infty$ and like $C  b^{-1/2}\rho^{-4}$ as $\rho\rightarrow\infty$. 

We will solve the differential equation $A\psi=0$ numerically in the region $(-\rho_0,\rho_0)$ and estimate the solution in the asymptotic region $|\rho|>\rho_0$. 

We proceed with the asymptotic analysis first, with the goal of finding conditions on the numerical solutions to guarantee that the solution will remain positive in the asymptotic region. Since $-gV'$ is positive for large $\rho$, we will ignore it in the asymptotic region, leaving this part available to be used as $\chispec$.  

To investigate asymptotic behaviour, we use the Schrodinger operator $B$ in terms of the variable $x$. We analyse the situation for $\rho\gg 1$, but the case $\rho\ll -1$ is similar. In this paragraph (and only this paragraph), we use $a \sim b$ to mean there is an $\epsilon$ such that $(1-\epsilon) a< b< (1+\epsilon) a$, with $\epsilon<2/(1+b\rho_0^2)$. 
\begin{align*}
\frac{d\rho}{dx} =& g'\sim \frac{1}{b\rho^2} \\
b\rho^2 \frac{d\rho}{dx} \sim& 1 \\
b\rho^3 \sim& 3x \\
g'''\sim& \frac{6}{b\rho^4} \\
-g'g'''/2
\sim& -\frac{6}{b^2\rho^6}\frac{1}{2} \\
\sim& -\frac{3}{9x^2} 
\end{align*}
Thus, taking $-(1+\epsilon)/3x^2$ as a lower bound for $-g'g'''/2$, it is sufficient to show that
\begin{align*}
B_1=-(2-\epsilon)\dx^2 - (1+\epsilon)\frac{1}{3x^2}
\end{align*}
has a positive solution to $B\psi =0$. Using $\epsilon'=3\epsilon<1/100$, we can replace this condition by positivity of the corresponding function for 
\begin{align*}
B=-2\dx^2 - \frac{1+\epsilon'}{3x^2} .
\end{align*}

Taking the ansatz $\phi=x^\alpha$, we can find solutions
\begin{align*}
\alpha=& \frac{ -2 \pm\sqrt{4-{8+8\epsilon'}/3}}{-4}= \frac12 \pm \frac{1}{2\sqrt{3}} +\text{$O(\epsilon')$ corrections} .
\end{align*}

A condition must be found to fit the asymptotics to the numerical solution. In the region $x\rightarrow\infty$, we require the positivity of the coefficient, $C_1$, on the more rapidly growing monomial, so that the solution will remain positive. Using
\begin{align*}
\phi(x)=& C_1 x^{\alpha_1 } + C_2 x^{\alpha_2 } \\
\phi'(x)=& C_1 \alpha_1 x^{\alpha_1-1 } + C_2 \alpha_2 x^{\alpha_2-1 } \\
C_2=& \frac{\phi(x) -C_1 x^{\alpha_1}}{x^{\alpha_2}} \\
\phi'(x)=& C_1 \alpha_1 x^{\alpha_1-1 } +\alpha_2\left(\phi(x)x^{-1} -C_1 x^{\alpha_1-1}\right) \\
C_1=& \frac{\phi'(x)- \alpha_2 \phi(x)x^{-1}}{\alpha_1-\alpha_2} x^{-\alpha_1+1} .
\end{align*}
Thus the condition we require is that, at the point $x$ where we match the numerics to the asymptotics, 
\begin{align*}
\dx \phi >&\alpha_2\phi(x) x^{-1} .
\end{align*}
Thus, we require(again with $\epsilon<2/(1+b\rho_0^2)$) 
\begin{align}
g'\dr \psi(\rho) >&\frac{1}{1-\epsilon} \alpha_2\psi(\rho)\frac{3}{b\rho^3}\nonumber\\
\dr \psi(\rho) >& \frac{1}{1-\epsilon} \alpha_2 \psi(\rho) \frac{3}{\rho} . \label{eMatchingCondition}
\end{align}
A similar result is required as $\rho\rightarrow-\infty$, but with the signs reversed. 

We now show numerically that there is a positive solution to
\begin{align*}
-(2-\epsilonspec)\dr g'\dr -g'''/2 -gV'
\end{align*}
on $(-\rho_0,\rho_0)$ with
\begin{align*}
\rho_0=&1000, &
\epsilonspec=& 1/1000,& 
b=& .1
\end{align*}
(we also take $M=1$ as a normalisation). This solution we construct satisfies 
\begin{align}
\dr \psi(\rho_0) \geq& 2(\frac12 -\frac{1}{2\sqrt3}) \psi(\rho_0) \frac{3}{\rho} , \nonumber\\
\dr \psi(-\rho_0) \geq& -2(\frac12 -\frac{1}{2\sqrt3}) \psi(-\rho_0) \frac{3}{\rho} , \label{eMatchAtLeft}
\end{align}
which, given the small size of $\epsilon=2/(1+b\rho^2)\sim 1/50,000$ and the additional factor of $2$ relative to \eqref{eMatchingCondition}, is sufficient to guarantee the matching of the numerics to the asymptotic solutions. 

We do this by treating $A\psi=0$ as an initial value problem for $\psi(\rho)$ posed at $\rho=\rho_0$ with initial conditions $\psi(\rho_0)=1$ and $\psi'(\rho_0)=2(1/2 -1/(2\sqrt3)) \psi(\rho_0) \frac{3}{\rho}$, finding the solution is positive, and verifying the condition is satisfied at $-\rho_0$. We show four plots: the potential $-g'''/2-gV'$, the solution in the left region $[0,\rho_0]$, the solution in the middle region $[-10,15]$ in which the potentials are large and the solution oscillates, and the solution in the right region $[-\rho_0,0]$ in which the solution goes rapidly to $\infty$ as $\rho\rightarrow-\infty$. Note that, although the asymptotic solution is concave down (with exponent less than 1) in terms of the variable $x$, since $x=b\rho^3/3$, the solution is concave up in terms of $\rho$ (and this behaviour is already clear from the plots as $\rho\rightarrow\pm\rho_0$). From the numerics, we find that at $-\rho_0=-1000$, the solution has value $\sim 150000$ and derivative $\sim -370$. From \eqref{eMatchAtLeft}, the derivative must be less than $\sim -170$. Since this is satisfied, the solution continued to the left of the numerical approximation will always be positive. By our choice of initial conditions, to the right of the numerical approximation, the continuation of the solution will also be positive. Thus, we have verified the ODE condition, condition \ref{cODE}, and hence the $L^2_{\text{loc}}$ result in \eqref{eMainResult}.

\begin{figure}
\includegraphics[width=3in]{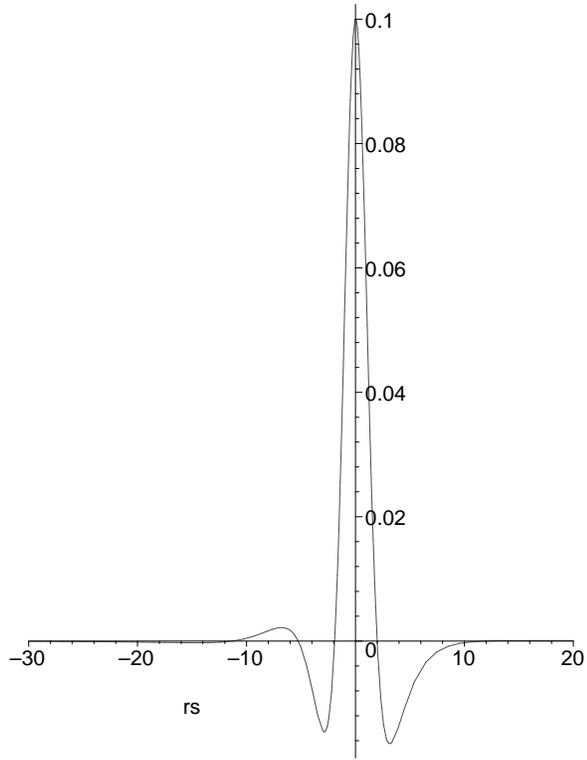}
\caption{The potentials $-g'''-gV'$.}
\end{figure}

\begin{figure}
\includegraphics[width=3in]{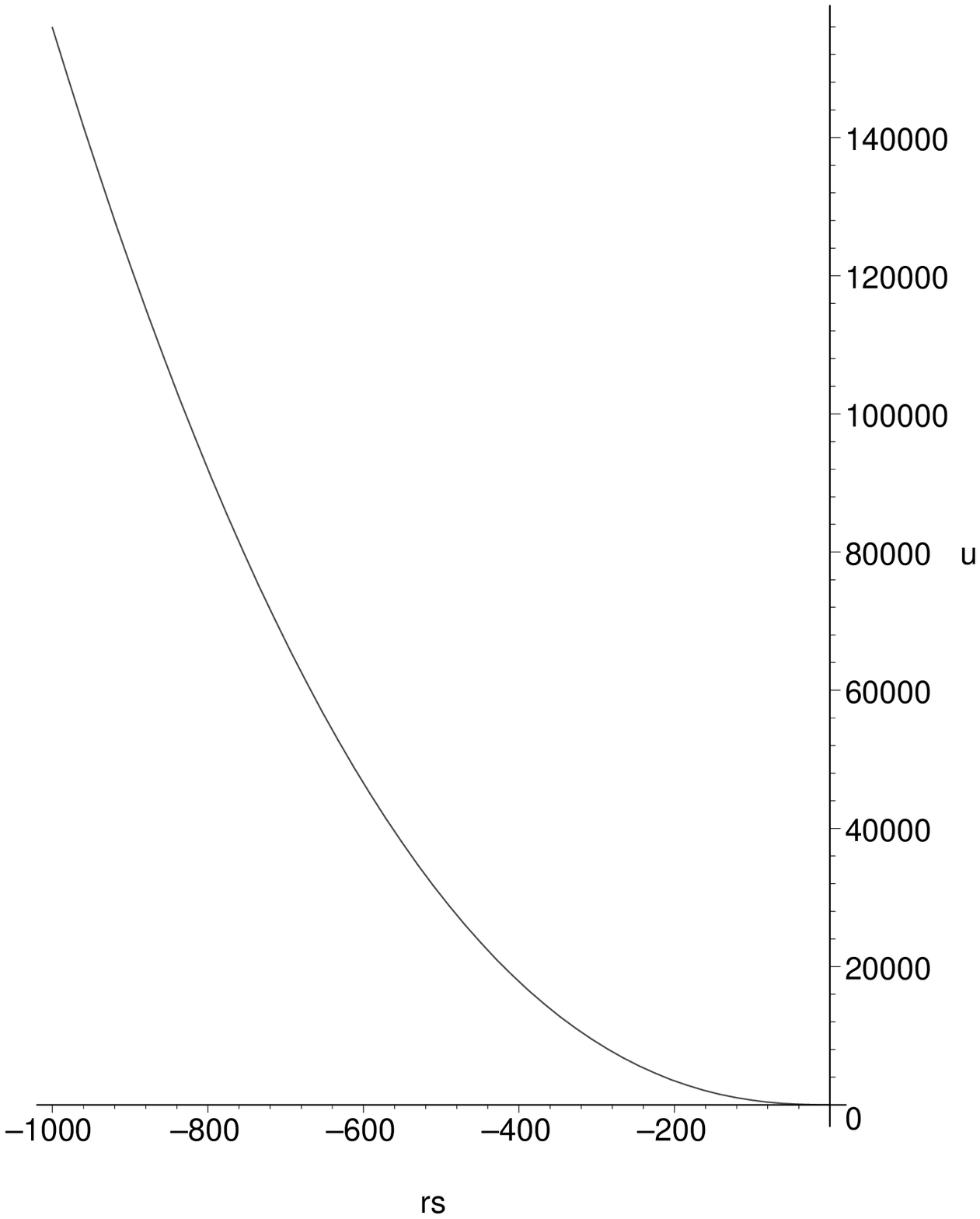}
\caption{The numerical solution in $[-1000,0]$.}
\end{figure}

\begin{figure}
\includegraphics[width=3in]{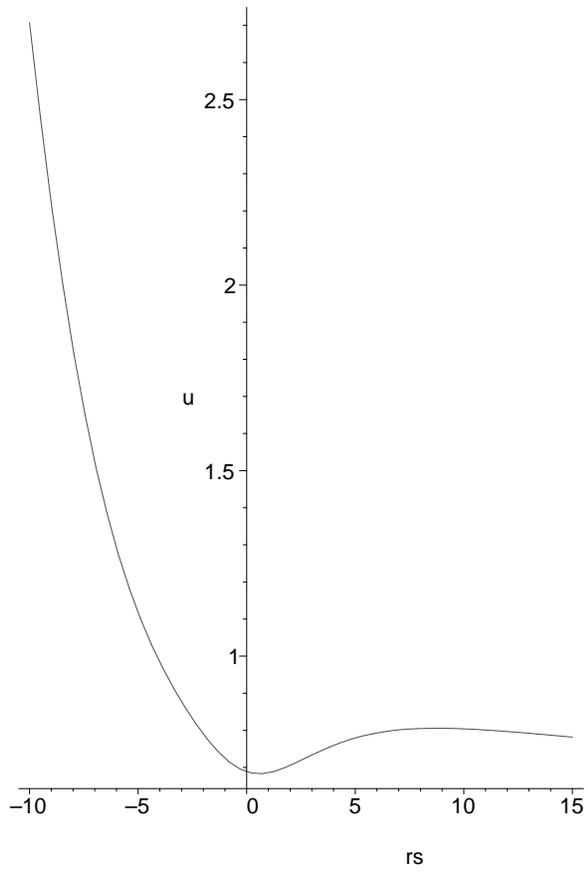}
\caption{The numerical solution in $[-10,15]$.}
\end{figure}

\begin{figure}
\includegraphics[width=3in]{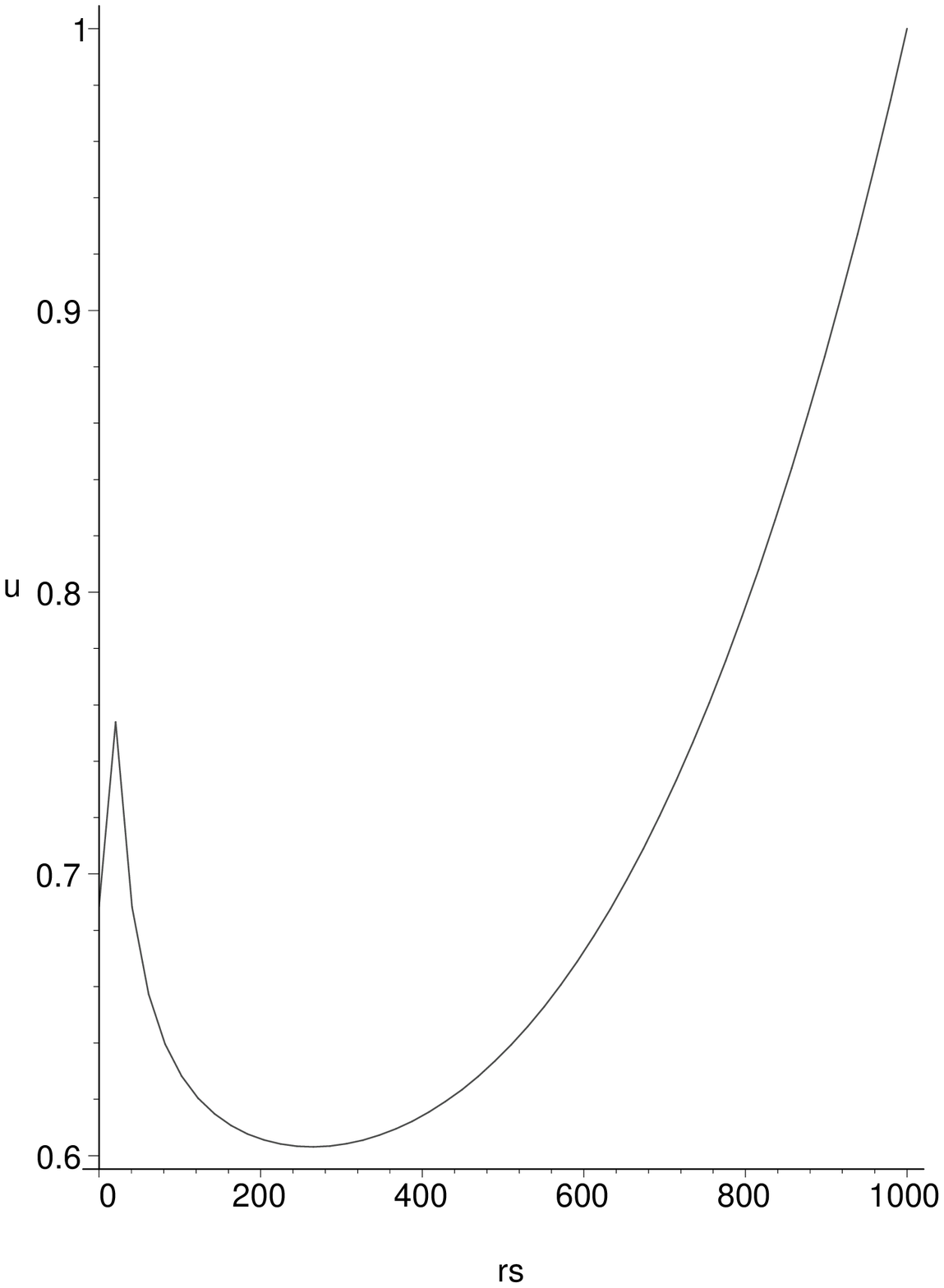}
\caption{The numerical solution in $[0,1000]$.}
\end{figure}


\vspace{.25in}
\noindent{\bf \large Acknowledgements}

\noindent A.S. is partially supported by NSF grant DMS-0501043. A.S. thanks Y. Pinchover for discussions. 



\begin{thebibliography}{1}
\bibitem{BlueSoffer3} P. Blue and A. Soffer: \emph{Phase Space Analysis on some Black Hole Manifolds} math.AP/0511281. 
\bibitem{BlueSoffer4} P. Blue and A. Soffer: \emph{Improved decay rates with small regularity loss for the wave equation about a Schwarzschild black hole.} math.AP/0612168. 
\bibitem{BlueSterbenz} P. Blue and J. Sterbenz: \emph{Uniform decay of local energy and the semi-linear wave equation on Schwarzschild space.}  Comm. Math. Phys. {\bf 268} (2006),  no. 2, 481--504.
\bibitem{DafermosRodnianski} M. Dafermos an I. Rodniansk:i \emph{The red-shift effect and radiation decay on black hole spacetimes.} gr-qc/0512119
\end{thebibliography}
\end{document}